\documentclass[12pt, reqno]{amsart}
\usepackage{amsfonts}
\usepackage{amsmath,amssymb,amscd,latexsym}
\usepackage{mathrsfs}
\usepackage{tikz-cd}

\newtheorem{theorem}{Theorem}[section]
\newtheorem{proposition}[theorem]{Proposition}

\newtheorem{lemma}[theorem]{Lemma}
\newtheorem{corollary}[theorem]{Corollary}
\newtheorem{remark}[theorem]{Remark}
\newtheorem{example}[theorem]{Example}

\begin{document}

\title{On motivic multiplicativity of cyclic covers}

\author{Ze Xu}

\date{October 31, 2026}

\begin{abstract}
Let $Y$ be a smooth projective complex variety of dimension $n$ with trivial Chow groups, and let $f: X\rightarrow Y$ be a cyclic cover of degree $d$ branched over a smooth ample divisor $D\in |dL|$.
We show that a naturally constructed Chow--Künneth decomposition of $X$ is $m$-fold $mn$-multiplicative for each integer $m\geq2$ and that its stable motivic multiplicativity defect satisfies $\mathfrak{sd}(X)\leq 2dn$.
If, in addition, each positive-codimension cycle class on $Y$ is supported on $D$, then a double cover admits motivic $0$-multiplicativity. As a further application, we show that any smooth Gushel--Mukai sixfold admits motivic 0-multiplicativity. Moreover, we establish the sharp modified-diagonal property for all smooth double EPW sextics, 
confirming a conjecture of O'Grady for these varieties.
\end{abstract}

\maketitle

\section{Introduction}
\bigskip

Let $X$ be a smooth conneced projective complex variety of dimension $n$, 
equipped with a Chow--K\"unneth (CK) decomposition $\{\pi_i^X\}_{0\leq i\leq 2n}$. 
With $\mathbb{Q}$-coefficients, denote
\begin{equation*}
\mathrm{CH}^r_s(X):=(\pi_{2r-s}^X)_*\mathrm{CH}^r(X),
\qquad
\pi_j^{X^{\ell}}:=\sum_{i_1+\cdots+i_{\ell}=j}\pi_{i_1}^X\otimes\cdots\otimes\pi_{i_{\ell}}^X.
\end{equation*}
As defined in \cite{X25}, this CK decomposition is $m$-fold $\tau$-multiplicative for $m\geq 2$ and $\tau\geq 0$, 
if $\Delta_{I_m}^X\in\bigoplus_{s=0}^{\tau}\mathrm{CH}^{mn}_s(X^{m+1})$.
Write $\mathfrak{d}_m(X)$ for the least such $\tau$ among all CK decompositions and $\mathfrak{sd}(X):=\sup_{m\geq2}\mathfrak{d}_m(X)$.
Motivic $0$-multiplicativity means the existence of a multiplicative CK decomposition.

Let $o_X$ be a degree-one 0-cycle class on $X$.
Put $I_k:=\{1,\ldots,k+1\}$ for $k\geq0$, and write $\delta_i^X:X\longrightarrow X^i$ for the $i$-th diagonal embedding for $i\geq 1$. Thus $\Delta_{I_k}^X:=(\delta_{k+1}^X)_*[X]$ and $\Delta_{I_1}^X=[X]$.
For each integer $k\geq3$, denote the $k$-th modified-diagonal class by
\begin{equation*}
\Gamma^k(X,o_X):=\sum_{J\subsetneq\{1,\ldots,k\}}(-1)^{|J|}p_J^*(o_X^{*|J|})\cdot p_{J'}^*\Delta_{I_{k-|J|-1}}^X
\in\mathrm{CH}_n(X^k),
\end{equation*}
where $J'$ is the complement of $J$, $p_J$ and $p_{J'}$ are the corresponding projections, and $o_X^{*j}$ is the exterior product of $j$ copies of $o_X$.

In \cite{X25}, it is conjectured that $\mathfrak{d}_m(X)\leq mn$ for each integer $m\geq 2$,
which has been verified for curves, surfaces, and ample subvarieties of varieties with trivial Chow groups.
In this companion paper, we investigate this conjecture and establish a bound on the stable motivic multiplicativity defect for cyclic covers of varieties with trivial Chow groups.

Our main result can be stated as follows.
\begin{theorem}\label{thm1.1}
Let $Y$ be a smooth connected projective variety of dimension $n$, and let $L$ be an ample divisor on $Y$. 
Let $f: X\longrightarrow Y$ be a cyclic cover of degree $d\geq 2$, branched along a smooth divisor $D\in |dL|$.
Suppose that $Y$ has trivial Chow groups.
Let $R\subseteq X$ be the reduced ramification divisor and let $o_X$ be the class of any point of $R$.
For $k\geq1$, set $f_k:=f^{\times k}$ and $R^*(X^k):=f_k^*\mathrm{CH}^*(Y^k)$.
\begin{enumerate}

\item[(i)] The natural self-dual CK decomposition constructed below is $m$-fold $mn$-multiplicative for each integer $m\geq2$. Moreover, $\mathfrak{sd}(X)\leq 2dn$.

\item[(ii)] Assume that each cycle class $z\in\mathrm{CH}^r(Y)$ of codimension $r\geq 1$ is supported on the branch divisor $D$. This holds, for example, for any smooth weighted projective space.
Then there exists an equality of cycle classes:
\begin{equation}\label{1}
\Gamma^{d+1}(X,o_X)+f_{d+1}^*w=0
\end{equation}
in $\mathrm{CH}^{dn}(X^{d+1})$,
where $f_{m}: X^m\rightarrow Y^m$ denotes the $m$-fold power of $f$ for each $m\geq 1$, and where $w\in\mathrm{CH}^{dn}(Y^{d+1})$. 
In particular, if $d=2$, then $X$ admits motivic 0-multiplicativity.

\item[(iii)] If $d=2$, the natural self-dual CK decomposition  is multiplicative if and only if
$\delta_{X*}f^*z\in R^*(X^2)$ for each cycle class $z\in\mathrm{CH}^r(Y)$ with $r\geq1$, 
where $\delta_X:=\delta_2^X$.
\end{enumerate}
\end{theorem}
 
As applications of further geometric arguments, we deduce the following results: 

(1) Let $X$ be a smooth projective surface with a relatively minimal elliptic fibration $X\longrightarrow\mathbb P^1$ admitting a section. Then $X$ admits motivic $0$-multiplicativity.

(2) Any smooth Gushel--Mukai sixfold admits motivic 0-multiplicativity. 
This improves a result in \cite{BL24}, which was only proved modulo algebraic equivalence.

(3) Any smooth double EPW sextic satisfies the sharp modified-diagonal property. 
This confirms a conjecture of O'Grady \cite{O'G14} for these hyper-Kähler varieties.
\bigskip

\textbf{Conventions}. We work over the field $\mathbb{C}$ of complex numbers unless otherwise stated.
The word \emph{variety} will refer to a reduced separated scheme of finite type over $\mathbb{C}$.
Singular cohomology groups, Chow groups, Chow rings and Chow motives are all with $\mathbb{Q}$-coefficients.
Denote by $\mathfrak{h}: \mathcal{V}(\mathbb{C})\rightarrow\mathrm{CH}\mathcal{M}$ the contravariant functor from the category of smooth complex projective varieties to the rigid pseudo-abelian tensor category of Chow motives with $\mathbb{Q}$-coefficients.
\bigskip

\section{Proof of Main theorem}
\bigskip

In this section, we prove Theorem \ref{thm1.1}. The key argument makes full use of combinatorics method encoded by the geometry of cyclic covers.

Write $b_{2i}:=\dim_{\mathbb Q}H^{2i}(Y,\mathbb Q)=\dim_{\mathbb Q}\mathrm{CH}^i(Y)$, and choose a basis $\{\beta_{ij}\}_{j=1}^{b_{2i}}$ of $\mathrm{CH}^i(Y)$ and its Poincar\'e dual basis $\{\beta_{ij}^{\vee}\}$ in $\mathrm{CH}^{n-i}(Y)$, so that $\deg(\beta_{ij}\beta_{ij'}^{\vee})=\delta_{jj'}$.
We always take $\beta_{01}=[Y]$ and $\beta_{n1}=o_Y$, the unique degree-one class in $\mathrm{CH}^n(Y)$.
Set
\begin{equation*}
\pi_{2i}^Y:=\sum_{j=1}^{b_{2i}}\beta_{ij}^{\vee}\times\beta_{ij},\qquad \pi_{2i+1}^Y:=0.
\end{equation*}
These are self-dual CK projectors of $Y$. Define
\begin{equation*}
\pi_i^X:=\frac1d f_2^*\pi_i^Y\quad(i\neq n),\qquad
\pi_n^X:=\Delta_X-\sum_{i\neq n}\pi_i^X.
\end{equation*}
Then $\{\pi_i^X\}_{0\leq i\leq 2n}$ forms a self-dual CK decomposition of $X$, by the generalized Lefschetz theorem \cite{F87}.
We will call it \emph{natural}.
Note that $f^*o_Y=do_X$ and each ambient cycle class $f^*z$ has pure grade-0.

The proof of Theorem 1.1 proceeds as follows.

\begin{proof}
(i) We first prove the bound $mn$ directly. Apply a product of the projectors to $\Delta_{I_m}^X$.
Observe that applying a nonmiddle projector replaces its factor by an ambient class and multiplies the class on the remaining diagonal by an ambient cycle class.
After applying all such projectors, the remaining term, if there remain $k\geq 1$ factors, has the form
$\bigl((\pi_n^X)^{\otimes k}\bigr)_*(\delta_k^X)_*\beta$, intersecting with grade-zero ambient classes, where $\beta\in f^*\mathrm{CH}^r(Y)$ and $0\leq r\leq n$; terms with $r>n$ vanish.
For $k\geq2$, its grade is $(k-2)n+2r$, which lies between $0$ and $kn$.
If at least one nonmiddle projector was applied, then $k\leq m$. If all $m+1$ factors are middle, then $r=0$ and the grade is $(m-1)n$.
For $k=1$, the remaining ambient class has grade zero or is killed by $\pi_n^X$; for $k=0$, a nonzero term is numerical and has grade zero.
Consequently, $\Delta_{I_m}^X$ has grades only in $[0,mn]$ and thus $\mathfrak d_m(X)\leq mn$.

Let $\sigma: X\rightarrow X$ be the covering involution. Clearly,
\begin{equation}\label{2}
f_2^*\Delta_Y=\Delta_X+\sum_{i=1}^{d-1}\Gamma_{\sigma^i}
\end{equation}
in $\mathrm{CH}^n(X\times X)$.
If $z\in\mathrm{CH}^r(Y)$ is supported on $D$, then for every $k\geq1$, there is a cycle equality
\begin{equation*}
(\delta_k^X)_*f^*z=d^{1-k}f_k^*(\delta_k^Y)_*z.
\end{equation*}
By hard Lefschetz and the assumption on $Y$, each cycle class of codimension $r>\frac{n}{2}$ lies in
$L^{2r-n}\mathrm{CH}^{n-r}(Y)$ and is supported on $D$.

To prove the stable bound, put
\begin{equation*}
\Pi:=\frac1d f_2^*\Delta_Y=\frac1d\sum_{j=0}^{d-1}\Gamma_{\sigma^j},\qquad
I:=(X,\Pi),\qquad W:=(X,\Delta_X-\Pi).
\end{equation*}
Observe that the invariant motive $I\cong\mathfrak h(Y)$ is a graded algebra with Tate summands $I^r$ of cohomological degree $2r$, and $W$ has cohomology only in degree $n$.
We may extend coefficients to $\mathbb Q(\zeta_d)$ and decompose $W$ into the nontrivial character summands $W_{\chi}$.
By the previous equality, we deduce that
\begin{enumerate}
\item[(a)] $I^r\otimes W\longrightarrow\mathfrak h(X)$ is zero for $r>\frac{n}{2}$;
\item[(b)] A product of $k\geq2$ nontrivial character summands with total character zero has no component in $I^r$ for $r<\frac{n}{2}$.
\end{enumerate}
For (a), the multiplication by each ambient cycle class of codimension $r$ is represented by its weighted two-factor diagonal, which is ambient and is killed by a $W$ projector.
For (b), dualize the projection to $I^r$: its components are the weighted $b$-factor diagonals of ambient classes of codimension $n-r>n/2$. These are ambient cycle class and are killed by the $W$ projectors.

Consider an arbitrary $m$-fold product of CK summands, with $s$ inputs in $W$ and invariant inputs of total codimension $A$.
If $s=0$, the product has defect zero. If $s>0$, (a) implies that a nonzero product has $A\leq\frac{n}{2}$.
A sequence of nonzero characters of a cyclic group of order $d$ can be partitioned into minimal zero-sum blocks of length at most $d$ and a zero-sum-free remainder of length at most $d-1$.
The length bounds follow by comparing partial sums of a minimal zero-sum sequence, or of a zero-sum-free sequence.
By equivariance and associativity, first multiply within each zero-sum block. By (b), its nonzero outputs have invariant codimension $r_j\geq\frac{n}{2}$.
If the final output is invariant, the remainder is empty, and $A+\sum_j r_j\leq n$ allows at most two blocks. Thus $s\leq2d$ and the motivic multiplicativity defect is precisely $sn-2\sum_jr_j\in[0,2dn]$.
Nonnegativity follows because every block has at least two inputs and $r_j\leq n$.
If the final output is noninvariant, two zero-sum blocks would give invariant codimension at least $n$ and annihilate the remaining $W$ input by (a).
Therefore, there is at most one such block and so $s\leq2d-1$. 
The defect is precisely $(s-1)n+2A\in[0,(2d-1)n]$.
These bounds hold for each $m\geq 2$. Hence, the natural CK decomposition is $m$-fold $2dn$-multiplicative for each $m\geq 2$, and $\mathfrak{sd}(X)\leq2dn$.

(ii) We next establish the cycle equality \eqref{1}.
Note that $\Delta_{I_k}^Y=\prod_{i=1}^{k}(p_{i,i+1}^Y)^*\Delta_Y$. 
For $1\leq k\leq d$ and $0\leq i_{\ell}\leq d-1$, define
\begin{equation*}
\Gamma_{i_0i_1\cdots i_k}:=\{(\sigma^{i_0}(x),\sigma^{i_1}(x),\cdots,\sigma^{i_k}(x))|x\in X\}\subseteq X^{k+1}.
\end{equation*}

Note that if $(i_1,\cdots, i_k)\neq (i'_1,\cdots,i'_k)$, then $\Gamma_{0i_1\cdots i_k}\neq\Gamma_{0i'_1\cdots i'_k}$.
For subsequent uses, observe that $\Gamma_{0i_1\cdots i_k}=\Gamma_{ii'_1\cdots i'_k}$ for each $0\leq i\leq d-1$, 
where for each $j$, $i'_j$ is the unique integer satisfying $0\leq i'_j\leq d-1$ and $i'_j\equiv i_j+i\pmod d$.
It follows that
\begin{equation}\label{3}
f_{k+1}^*\Delta_{I_k}^Y
=\prod_{r=1}^{k}(p_{r,r+1}^X)^*\left(\Delta_X+\sum_{j=1}^{d-1}\Gamma_{\sigma^j}\right)
=\sum_{i_1=0}^{d-1}\cdots\sum_{i_k=0}^{d-1}\Gamma_{0i_1\cdots i_k}.
\end{equation}

Let $\mathcal{P}_{k+1}$ denote the set of partitions of $k+1$, and let $\mathscr{P}_{k+1}$ denote the set of set partitions of $\{1, 2,\cdots, k+1\}$. 
For $\mathbf{p}\in\mathscr{P}_{k+1}$, let $|\mathbf{p}|$ denote the number of blocks of $\mathbf{p}$.
For any partition $\mathfrak{a}=1^{a_1}2^{a_2}\cdots(k+1)^{a_{k+1}}\in\mathcal{P}_{k+1}$ with each $a_{\ell}\geq 0$,
let $M_{\mathfrak{a}}$ denote the set of multi-sets $\{0,i_1,\cdots, i_k\}$ of integers such that each $0\leq i_j\leq d-1$ and there are exactly $a_{\ell}$ elements in its underlying set,
denoted $\{0,i_1,\cdots, i_k\}_{\text{red}}=\{i'_1,\cdots, i'_{l_{\mathfrak{a}}}\}$, each with multiplicity $\ell$ for each $\ell$, 
where $l_{\mathfrak{a}}:=\sum_{\ell=1}^{k+1}a_{\ell}$ is the length of $\mathfrak{a}$. 
Define 
\begin{equation*}
F_{\mathfrak{a}}:=\{\mathbf{p}\in\mathscr{P}_{k+1}|\mathbf{p}\ \text{has exactly}\ a_{\ell}\ \text{blocks of size}\ \ell\},
\end{equation*}
\begin{equation*}
S_{\mathfrak{a}}:=\{(i_1,\cdots, i_k)\in\mathbb{Z}^k|\{0,i_1,\cdots, i_k\}\in M_{\mathfrak{a}}\}
\end{equation*}
and
\begin{equation*}
\widehat{S}_{\mathfrak{a}}:=\{(0,i'_2,\ldots,i'_{l_{\mathfrak{a}}})\in\mathbb Z^{l_{\mathfrak{a}}}
\mid \{i'_1,\ldots,i'_{l_{\mathfrak{a}}}\}
 =\{0,i_1,\ldots,i_k\}_{\mathrm{red}}\text{ for some }\{0,i_1,\ldots,i_k\}\in M_{\mathfrak{a}}\}.
\end{equation*}
Every $\mathbf{p}\in\mathscr{P}_{k+1}$ is naturally ordered by arranging its blocks according to their least elements.
Define a surjective map
\begin{equation*}
\phi: S_{\mathfrak{a}}\rightarrow \widehat{S}_{\mathfrak{a}},\quad
(i_1,\cdots, i_k)\mapsto (i'_1=0, i'_2,\cdots, i'_{l_{\mathfrak{a}}}),
\end{equation*}
where $i'_j$ is the $j$-th distinct entry of the vector $(0,i_1,\cdots, i_k)$.
For any $(i'_1,\cdots, i'_{l_{\mathfrak{a}}})\in\widehat{S}_{\mathfrak{a}}$, 
denote $F_{(i'_1,\cdots, i'_{l_{\mathfrak{a}}})}:=\phi^{-1}(i'_1,\cdots, i'_{l_{\mathfrak{a}}})$.
Define a map
\begin{equation*}
\psi: F_{(i'_1,\cdots, i'_{l_{\mathfrak{a}}})}\rightarrow F_{\mathfrak{a}},\quad
(i_1,\cdots, i_k)\mapsto\mathbf{p}_{(i_1,\cdots, i_k)},
\end{equation*}
where $\mathbf{p}_{(i_1,\cdots, i_k)}$ is the partition of the set of ordinal numbers of the vector $(0, i_1,\cdots, i_k)$ uniquely determined by the same components.
Obviously, $\psi$ is a bijection.
Thus, one knows that
\begin{equation*}
|F_{(i'_1,\cdots, i'_{l_{\mathfrak{a}}})}|
=|F_{\mathfrak{a}}|
=\frac{(k+1)!}{\prod_{\ell}a_{\ell}!(\ell!)^{a_{\ell}}}.
\end{equation*}
Recall that $k+1=\sum_{\ell=1}^{k+1}\ell a_{\ell}$.
Then the cardinalities of $\widehat{S}_{\mathfrak{a}}$, $M_{\mathfrak{a}}$ and $S_{\mathfrak{a}}$ are respectively
\begin{equation*}
|\widehat{S}_{\mathfrak{a}}|=\prod_{i=1}^{l_{\mathfrak{a}}-1}(d-i),\quad
|M_{\mathfrak{a}}|=\frac{l_{\mathfrak{a}}!}{\prod_{\ell}a_{\ell}!}\binom{d-1}{l_{\mathfrak{a}}-1},
\end{equation*}
\begin{equation*}
\begin{aligned}
|S_{\mathfrak{a}}|
&=|F_{\mathfrak{a}}||\widehat{S}_{\mathfrak{a}}| \\
&=\frac{(k+1)!}{\prod_{\ell}a_{\ell}!(\ell!)^{a_{\ell}}}\prod_{i=1}^{l_{\mathfrak{a}}-1}(d-i).
\end{aligned}
\end{equation*}
In particular, $|S_{(k+1)^1}|=1$ and $|S_{1^{k+1}}|=|\widehat{S}_{1^{k+1}}|=\prod_{i=1}^{k}(d-i)$.
It follows that
\begin{equation}\label{4}
d^k=\sum_{\mathfrak{a}\in\mathcal{P}_{k+1}}|S_{\mathfrak{a}}|
=\sum_{\mathfrak{a}\in\mathcal{P}_{k+1}}\frac{(k+1)!}{\prod_{\ell}a_{\ell}!(\ell!)^{a_{\ell}}}
\prod_{i=1}^{l_{\mathfrak{a}}-1}(d-i).
\end{equation}

For any $\mathfrak{a}\in\mathcal{P}_{k+1}$, define
\begin{equation*}
\Theta_{\mathfrak{a}}:=\sum_{(i_1,\cdots,i_k)\in S_{\mathfrak{a}}}\Gamma_{0i_1\cdots i_k}.
\end{equation*}
In particular,
\begin{equation*}
\Theta_{(k+1)^1}=\Delta_{I_k}^X,\quad
\Theta_{1^{k+1}}=\sum_{\text{distinct}\ i_j\geq 1}\Gamma_{0i_1\cdots i_k}\in\text{CH}^{kn}(X^{k+1}).
\end{equation*}
By \eqref{2}, $\Theta_{1^2}=f_2^*\Delta_Y-\Delta_X$, whereas $\Theta_{1^{d+1}}=0$.
For $\mathbf{p}\in F_{\mathfrak{a}}$, define the morphism
\begin{equation*}
\delta_{\mathbf{p}}^{X^{l_{\mathfrak{a}}}}: X^{l_{\mathfrak{a}}}\rightarrow X^{k+1},\quad
(x_1, \cdots, x_{l_{\mathfrak{a}}})\mapsto 
(x'_j)_{j=0}^k,
\end{equation*}
where $x'_j=x_{\ell}$ whenever $j+1$ belongs to the $\ell$-th block of $\mathbf{p}$.
Equivalently, for every $(i'_1=0,\cdots, i'_{l_{\mathfrak{a}}})\in\widehat{S}_{\mathfrak{a}}$ and every $\mathbf{p}\in F_{\mathfrak{a}}$, 
there is a unique tuple $(i_1,\cdots,i_k)\in F_{(i'_1,\cdots, i'_{l_{\mathfrak{a}}})}$ whose equality pattern is described by $\mathbf p$, and whose distinct entries, ordered by their first occurrence, are $i'_1=0,\cdots, i'_{l_{\mathfrak{a}}}$.
For every $\mathfrak{a}\in\mathcal{P}_{k+1}$, we therefore have
\begin{equation}\label{5}
\begin{split}
 \Theta_{\mathfrak{a}} &=\sum_{(i_1,\cdots,i_k)\in S_{\mathfrak{a}}}\Gamma_{0i_1\cdots i_k} \\
     &=\sum_{\mathbf{p}\in F_{\mathfrak{a}}}\sum_{(0,i'_2,\cdots,i'_{l_{\mathfrak{a}}})\in \widehat{S}_{\mathfrak{a}}}(\delta_{\mathbf{p}}^{X^{l_{\mathfrak{a}}}})_*\Gamma_{0i'_2\cdots i'_{l_{\mathfrak{a}}}}\\
     &=\sum_{\mathbf{p}\in F_{\mathfrak{a}}}(\delta_{\mathbf{p}}^{X^{l_{\mathfrak{a}}}})_*\Theta_{1^{l_{\mathfrak{a}}}}.
\end{split}
\end{equation}
Thus for each $1\leq k\leq d$, it follows from \eqref{3} and \eqref{5} that
\begin{equation}\label{6}
\begin{split}
   f_{k+1}^*\Delta_{I_k}^Y & =\sum_{\mathfrak{a}\in\mathcal{P}_{k+1}}\sum_{(i_1,\cdots,i_k)\in S_{\mathfrak{a}}}\Gamma_{0i_1\cdots i_k} \\
     &=\sum_{\mathfrak{a}\in\mathcal{P}_{k+1}}\Theta_{\mathfrak{a}}\\
     &=\sum_{\mathfrak{a}\in\mathcal{P}_{k+1}}\sum_{\mathbf{p}\in F_{\mathfrak{a}}}(\delta_{\mathbf{p}}^{X^{l_{\mathfrak{a}}}})_*\Theta_{1^{l_{\mathfrak{a}}}}.
\end{split}
\end{equation}
Hence,
\begin{equation}\label{7}
\Theta_{1^{k+1}}=f_{k+1}^*\Delta_{I_k}^Y-\sum_{\mathfrak{a}\in\mathcal{P}_{k+1}-\{1^{k+1}\}}\sum_{\mathbf{p}\in F_{\mathfrak{a}}}(\delta_{\mathbf{p}}^{X^{l_{\mathfrak{a}}}})_*\Theta_{1^{l_{\mathfrak{a}}}}.
\end{equation}
Equivalently, M\"obius inversion on the lattice of set partitions gives the uniform identity
\begin{equation*}
\Theta_{1^{k+1}}=\sum_{\mathbf p\in\mathscr P_{k+1}}\mu(\widehat0,\mathbf p)
(\delta_{\mathbf p}^{X^{|\mathbf p|}})_*f_{|\mathbf p|}^*\Delta_{I_{|\mathbf p|-1}}^Y,
\qquad
\mu(\widehat0,\mathbf p)=\prod_{B\in\mathbf p}(-1)^{|B|-1}(|B|-1)!.
\end{equation*}
Every term here is a class on $X^{k+1}$, including the one-block term, for which $f_1=f$ and $\Delta_{I_0}^Y=[Y]$.
For $0\leq i\leq k$, put
\begin{equation*}
S_{k,i}:=\sum_{\substack{J\subset I_k\\ |J|=i}}p_{J'}^*\Delta_{I_{k-i}}^X\cdot p_J^*(o_X^{*i}).
\end{equation*}
Thus each choice of the point factors occurs once.
Expanding the diagonals of $Y$ into exterior products and separating terms with a fundamental-class factor yields
\begin{equation}\label{8}
\Theta_{1^{k+1}}=\sum_{i=0}^{k}c_{ki}S_{k,i}
+\sum_{\substack{\mathfrak a\in\mathcal P_{k+1}\\ l_{\mathfrak a}\geq2}}
\sum_{\mathbf p\in F_{\mathfrak a}}(\delta_{\mathbf p}^{X^{l_{\mathfrak a}}})_*f_{l_{\mathfrak a}}^*
\Phi_{\mathfrak a,\mathbf p}^{l_{\mathfrak a}}.
\end{equation}
Here $\Phi_{\mathfrak a,\mathbf p}^{l_{\mathfrak a}}\in\mathrm{CH}^{(l_{\mathfrak a}-1)n}(Y^{l_{\mathfrak a}})$ is a $\mathbb{Q}$-linear combination of exterior products whose factors all have positive codimension, allowing codimension $n$.
Indeed, a decomposable term of codimension $(l-1)n$ with a fundamental-class factor has all its other factors in codimension $n$.
Its pullback and diagonal pushforward contribute to the $S_{k,i}$, because $f^*o_Y=do_X$ and $(\delta_b^X)_*o_X=o_X^{*b}$.
The remaining terms have positive codimension on every block; their further diagonal pushforwards still have the displayed form on partitions of $I_k$.
Symmetry makes $c_{ki}$ independent of the particular $i$ point factors.
The number of partitions $\mathbf{p}\in F_{\mathfrak{a}}$ for which exactly $j$ blocks have union $\{1,2,\cdots, i\}$ is precisely
\begin{equation*}
\frac{i!(k+1-i)!}{\prod_{\ell}a_{\ell}!(\ell!)^{a_{\ell}}}\cdot A_{\mathfrak{a}}(j,i)
=\sum_{\substack{0\leq b_{\ell}\leq a_{\ell}\\
\sum_{\ell}b_{\ell}=j\\
\sum_{\ell}\ell b_{\ell}=i}}
\frac{i!}{\prod_{\ell}b_{\ell}!(\ell!)^{b_{\ell}}}\cdot\frac{(k+1-i)!}{\prod_{\ell}(a_{\ell}-b_{\ell})!(\ell!)^{a_{\ell}-b_{\ell}}},
\end{equation*}
where 
\begin{equation*}
A_{\mathfrak{a}}(j,i):=\sum_{\substack{0\leq b_{\ell}\leq a_{\ell}\\
\sum_{\ell}b_{\ell}=j\\
\sum_{\ell}\ell b_{\ell}=i}}
\prod_{\ell}\binom{a_{\ell}}{b_{\ell}}
\end{equation*}
is the coefficient of $x^jy^i$ in the polynomial $P_{\mathfrak{a}}(x,y):=\prod_{\ell}(1+xy^{\ell})^{a_{\ell}}$.
Now for any integer $k\geq 0$ and $i\geq 0$, set $c_{00}=1$ and $c_{ki}=0$ for $i>k$. 
Observe that $c_{10}=-1$, $c_{11}=d$ and there are the following recurrence relations:
\begin{equation}\label{9}
c_{ki}=-\sum_{\mathfrak{a}\in\mathcal{P}_{k+1}-\{1^{k+1}\}}
\sum_{j=0}^{i}c_{l_{\mathfrak{a}}-1,j}\frac{i!(k+1-i)!}
{\prod_{\ell}a_{\ell}!(\ell!)^{a_{\ell}}}\cdot A_{\mathfrak{a}}(j,i)
\end{equation}
for $k\geq 2$ and $i\leq k-1$, and
\begin{equation}\label{10}
c_{kk}=d^k-\sum_{\mathfrak{a}\in\mathcal{P}_{k+1}-\{1^{k+1}\}}
c_{l_{\mathfrak{a}}-1,l_{\mathfrak{a}}-1}\frac{a_1k!}{\prod_{\ell}a_{\ell}!(\ell!)^{a_{\ell}}}
\end{equation}
for $k\geq 2$, since $A_{\mathfrak{a}}(j, k)=\delta_{j,l_{\mathfrak{a}}-1}a_1$. 
Thus Lemma \ref{lem2.3} gives  
\begin{equation*}
c_{ki}=(-1)^{k-i}\prod_{j=0}^{i-1}(d-j)(k-i)!
\end{equation*}
for every $0\leq i\leq k$, with the empty product equal to $1$.
In particular, $c_{di}=(-1)^{d-i}d!$ for $0\leq i\leq d$.
Since $\Theta_{1^{d+1}}=0$, equation \eqref{8}, divided by $(-1)^dd!$, gives
\begin{equation}\label{11}
\Gamma^{d+1}(X,o_X)
+\frac{(-1)^d}{d!}\sum_{\substack{\mathfrak a\in\mathcal P_{d+1}\\ l_{\mathfrak a}\geq2}}
\sum_{\mathbf p\in F_{\mathfrak a}}(\delta_{\mathbf p}^{X^{l_{\mathfrak a}}})_*f_{l_{\mathfrak a}}^*
\Phi{}_{\mathfrak a,\mathbf p}^{l_{\mathfrak a}}=0.
\end{equation}
By hypothesis, every factor of each $\Phi$ is supported on $D$.
Applying the cycle equality proved in (i) separately to each block of $\mathbf p$, every correction term of \eqref{11} is an ambient cycle class on $X^{d+1}$.
Their sum is therefore $f_{d+1}^*w$ for some $w\in\mathrm{CH}^{dn}(Y^{d+1})$, proving \eqref{1}.
If $d=2$, the same cycle equality gives $\delta_{X*}f^*z\in R^*(X^2)$ for every positive-codimension $z$, and (iii) proves motivic $0$-multiplicativity.

(iii) For a double cover, $W$ is the anti-invariant summand, placed in degree $n$.
Assume first that every positive-codimension weighted diagonal is ambient.
Then multiplication $I^r\otimes W\longrightarrow W$ vanishes for $r>0$, since its correspondence is killed by the $W$ projector.
By equivariance, $W\otimes W$ maps to $I$. Its projection to $I^r$ for $r<n$ is dual to a weighted diagonal of codimension $n-r>0$ and is therefore zero.
Thus, $W\otimes W$ maps only to $I^n$, of degree $2n$. The multiplication on $I$ is graded and $I^0$ acts by the unit, so the natural CK decomposition is multiplicative.

Conversely, suppose that this CK decomposition is multiplicative.
For $z\in\mathrm{CH}^r(Y)$ with $r>0$, multiplication by $f^*z$ has no component from $W$ to $W$, since the input degree is $n+2r$, whereas the output degree would be $n$.
The weighted diagonal $\delta_{X*}f^*z$ is invariant under $\sigma\times\sigma$; its mixed $I\otimes W$ and $W\otimes I$ components vanish.
Its $W\otimes W$ component is precisely the multiplication component just excluded, by self-duality.
Therefore, it lies in the $I\otimes I$ summand, whose Chow groups are $R^*(X^2)$. This completes the proof.
\end{proof}

\begin{remark}\label{rem2.2}
\leavevmode
\begin{enumerate}
\item[(i)] For $z\in\mathrm{CH}^r(Y)$ with $r\geq1$, the following conditions are equivalent:
\begin{enumerate}
\item[(1)] $\delta_{X*}f^*z\in R^*(X^2)$.
\item[(2)] $\delta_{X*}f^*z=\frac1d f_2^*\delta_{Y*}z$.
\item[(3)] $\Gamma_\sigma\cdot(X\times f^*z)=\Delta_X\cdot(X\times f^*z)$.
\end{enumerate}

\item[(ii)] The natural CK construction, the bound $\mathfrak d_m(X)\leq mn$, and parts (ii) and (iii) remain valid for a finite flat cyclic cover of smooth projective varieties when the generalized Lefschetz property replaces ampleness: $f^*$ is an isomorphism outside the middle degree onto the cohomology of $X$.
One must also choose a point class $o_X$ satisfying $f^*o_Y=do_X$; this condition is automatic for a totally ramified point in the situation of the theorem.
The stable bound additionally uses that every class on $Y$ of codimension $>\frac{n}{2}$ is supported on the branch divisor. Ampleness supplies this property in Theorem \ref{thm1.1}.

\item[(iii)] Theorem \ref{thm1.1} also remains valid, when $Y$ and $D$ are merely normal varieties with at worst quotient singularities. This would resort to intersection theory of global quotients, refer to \cite{V89}.
\end{enumerate}
\end{remark}

\begin{lemma}\label{lem2.3}
For any integer $k\geq 0$ and $i\geq 0$, set $c_{00}=1$, $c_{10}=-1$, $c_{11}=d$ and $c_{ki}=0$ for $i>k$. 
Assume that for $k\geq 2$ and $i\leq k-1$, the integers $c_{ki}$ satisfy the recurrence relations \eqref{9} and \eqref{10}.
Then one has 
\begin{equation*}
c_{ki}=(-1)^{k-i}\prod_{j=0}^{i-1}(d-j)(k-i)!
\end{equation*}
for any $0\leq i\leq k$, with the empty product equal to $1$.
\end{lemma}
\begin{proof}
We proceed by induction on $k\geq 2$.
It is clear that $c_{20}=2$, $c_{21}=-d$ and $c_{22}=d(d-1)$.
To warm up, for $i=0$, one gets a recurrence relation:
\begin{equation*}
c_{k0}=-\sum_{\mathfrak{a}\in\mathcal{P}_{k+1}-\{1^{k+1}\}}c_{l_{\mathfrak{a}}-1,0}\frac{(k+1)!}
{\prod_{\ell}a_{\ell}!(\ell!)^{a_{\ell}}},
\end{equation*}
since $A_{\mathfrak{a}}(0, 0)=1$.
In the sums below, a term with a negative factorial argument is understood to be zero.
It follows that $c_{k0}=(-1)^kk!$, because by induction on $k$, one deduces that
\begin{eqnarray*}
 c_{k0} &=& \sum_{\mathfrak{a}\in\mathcal{P}_{k+1}-\{1^{k+1}\}}(-1)^{l_{\mathfrak{a}}}
 \frac{(k+1)!(l_{\mathfrak{a}}-1)!}{\prod_{\ell}a_{\ell}!(\ell!)^{a_{\ell}}}\\
   &=& \sum_{\mathbf{p}\in\mathscr{P}_{k+1}-\{\text{singletons}\}}(-1)^{|\mathbf{p}|}(|\mathbf{p}|-1)!\\
   &=& (-1)^kk!.
\end{eqnarray*}
The last equality makes use of the Möbius function formula, referring to page 344 in \cite{R64}.
By induction on $k\geq 2$ and the equality \eqref{4}, one deduces that 
\begin{eqnarray*}
c_{kk} &=& d^k-\sum_{\mathfrak{a}\in\mathcal{P}_{k+1}-\{1^{k+1}\}}
\frac{a_1k!}{\prod_{\ell}a_{\ell}!(\ell!)^{a_{\ell}}}\prod_{i=0}^{l_{\mathfrak{a}}-2}(d-i) \\
   &=&  \prod_{i=0}^{k-1}(d-i).
\end{eqnarray*}
For $i\leq k-1$, by induction on $k$, one gets that
\begin{eqnarray*}
 c_{ki}&=& -\sum_{\mathfrak{a}\in\mathcal{P}_{k+1}-\{1^{k+1}\}}
\sum_{j=0}^{i}c_{l_{\mathfrak{a}}-1,j}\frac{i!(k+1-i)!}
{\prod_{\ell}a_{\ell}!(\ell!)^{a_{\ell}}}\cdot A_{\mathfrak{a}}(j,i) \\
   &=& \sum_{\mathfrak{a}\in\mathcal{P}_{k+1}-\{1^{k+1}\}}
\sum_{j=0}^{i}(-1)^{l_{\mathfrak{a}}-j}\prod_{\ell=0}^{j-1}(d-\ell)(l_{\mathfrak{a}}-1-j)!\frac{i!(k+1-i)!}
{\prod_{\ell}a_{\ell}!(\ell!)^{a_{\ell}}}\cdot A_{\mathfrak{a}}(j,i).
\end{eqnarray*}
Note that $A_{1^{k+1}}(j,i)=\binom{k+1}{j}\delta_{ij}$ and
\begin{eqnarray*}
   &  & \sum_{j=0}^{i}(-1)^{k+1-j}\prod_{\ell=0}^{j-1}(d-\ell)(k-j)!\frac{i!(k+1-i)!}
{(k+1)!}\cdot A_{1^{k+1}}(j,i)\\
   &=& (-1)^{k+1-i}\prod_{\ell=0}^{i-1}(d-\ell)(k-i)!\frac{i!(k+1-i)!}
{(k+1)!}\binom{k+1}{i} \\
   &=& (-1)^{k+1-i}\prod_{\ell=0}^{i-1}(d-\ell)(k-i)!.
\end{eqnarray*}
Then it suffices to show that 
\begin{equation*}
\sum_{j=0}^{i}i!(k+1-i)!\prod_{\ell=0}^{j-1}(d-\ell)\cdot B(j,i)=0,
\end{equation*}
where 
\begin{equation*}
B(j,i):=\sum_{\mathfrak{a}\in\mathcal{P}_{k+1}}
(-1)^{l_{\mathfrak{a}}-j}(l_{\mathfrak{a}}-1-j)!\frac{A_{\mathfrak{a}}(j,i)}
{\prod_{\ell}a_{\ell}!(\ell!)^{a_{\ell}}}
\end{equation*}
for fixed $j$ and $i$. Then one is reduced to show that $B(j,i)=0$ for $i\leq k-1$.
For $\mathfrak{b}=1^{b_1}2^{b_2}\cdots\in\mathcal{P}_i$ such that $0\leq b_{\ell}\leq a_{\ell}$ and $l_{\mathfrak{b}}=j$,
denote $c_{\ell}:=a_{\ell}-b_{\ell}$ and thus $\mathfrak{c}:=1^{c_1}2^{c_2}\cdots\in\mathcal{P}_{k+1-i}$. 
Then
\begin{equation*}
B(j,i)=\sum_{\substack{\mathfrak{b}\in\mathcal{P}_i\\ l_{\mathfrak b}=j}}\prod_{\ell}\frac{1}{b_{\ell}!(\ell!)^{b_{\ell}}}\cdot
\sum_{\mathfrak{c}\in\mathcal{P}_{k+1-i}}(-1)^{l_{\mathfrak{c}}}(l_{\mathfrak{c}}-1)!
\prod_{\ell}\frac{1}{c_{\ell}!(\ell!)^{c_{\ell}}}.
\end{equation*}
It is enough to show that 
\begin{equation*}
C_{k+1-i}:=\sum_{\mathfrak{c}\in\mathcal{P}_{k+1-i}}
(-1)^{l_{\mathfrak{c}}}(l_{\mathfrak{c}}-1)!\prod_{\ell}\frac{1}{c_{\ell}!(\ell!)^{c_{\ell}}}=0.
\end{equation*}
Considering the corresponding generating function, one immediately calculates that
\begin{equation*}
G(t):=\sum_{i\geq 1}C_it^i=\sum_{i\geq 1}\sum_{\mathfrak{c}\in\mathcal{P}_{i}}
(-1)^{l_{\mathfrak{c}}}(l_{\mathfrak{c}}-1)!\prod_{\ell}\frac{1}{c_{\ell}!(\ell!)^{c_{\ell}}}t^i
=\sum_{i=1}^{\infty}(-1)^i\frac{(e^t-1)^i}{i}=-t.
\end{equation*}
Therefore, $C_{k+1-i}=0$ for $i\leq k-1$.
\end{proof}

As the first consequences of Theorem \ref{thm1.1}, we have the following.
\begin{example}
Let $Y=\mathbb P_w^n$ be a smooth weighted projective space, and let $f: X\longrightarrow Y$ be a smooth double cover branched along a smooth divisor $D\in|2L|$, where $L$ is ample and Cartier.
Since $\mathrm{CH}^*(Y)=\mathbb Q[H]/(H^{n+1})$ and $D$ is a nonzero multiple of $H$, all positive-codimension classes are supported on $D$.
Thus, it follows from Theorem \ref{thm1.1}(ii) that $X$ admits motivic $0$-multiplicativity.
\end{example}

\begin{example}
As an immediate consequence of Theorem \ref{thm1.1}, the Cynk--Hulek Calabi--Yau varieties constructed from elliptic curves using the order-two or order-three actions admit MCK decompositions, 
which were first proved in \cite{LV20} by using a different method.
\end{example}
\bigskip

\section{Further applications}
\bigskip
In this section, we draw several further applications of Theorem \ref{thm1.1}.

\bigskip

First, we specialize to surfaces equipped with involutions.

\begin{lemma}\label{lem3.1}
Let $X$ be a smooth projective regular surface, and let $C\subseteq X$ be an integral, possibly singular, constant-cycle curve.
Let $o_X\in\mathrm{CH}^2(X)$ be the degree-one class of any point of $C$.
Then $\delta_{X*}C=C\times o_X+o_X\times C$.
\end{lemma}
\begin{proof}
Let $\phi:\widetilde C\longrightarrow X$ be the normalization of $C$ followed by its inclusion, and choose a point class $o_{\widetilde C}$.
Set $\Gamma:=\Gamma_\phi\in\mathrm{CH}^2(\widetilde C\times X)$.
Use the self-dual projectors
\begin{equation*}
\pi_0^{\widetilde C}=o_{\widetilde C}\times\widetilde C,\quad
\pi_2^{\widetilde C}=\widetilde C\times o_{\widetilde C},\quad
\pi_1^{\widetilde C}=\Delta_{\widetilde C}-\pi_0^{\widetilde C}-\pi_2^{\widetilde C},
\end{equation*}
and
\begin{equation*}
\pi_0^X=o_X\times X,\quad \pi_1^X=\pi_3^X=0,\quad
\pi_4^X=X\times o_X,\quad \pi_2^X=\Delta_X-\pi_0^X-\pi_4^X.
\end{equation*}
A direct application of the product projectors, using regularity and $\phi_*o_{\widetilde C}=o_X$, gives
\begin{equation*}
\Gamma=\Gamma_0+\Gamma_1,\qquad
\Gamma_0=\widetilde C\times o_X+o_{\widetilde C}\times C,\qquad
\Gamma_1=(\pi_1^{\widetilde C}\otimes\pi_2^X)_*\Gamma.
\end{equation*}
Here $\Gamma_i\in\mathrm{CH}^2_i(\widetilde C\times X)$.
For every closed point $x\in\widetilde C$, the constant-cycle property gives
\begin{equation*}
\Gamma_{1*}x=(\pi_2^X)_*\phi_*(x-o_{\widetilde C})=0.
\end{equation*}
The decomposition principle of \cite[Theorem~10.19]{Voi03}, followed by localization, therefore gives $\Gamma_1=\sum_j x_j\times D_j$, with points $x_j\in\widetilde C$ and divisors $D_j$ on $X$.
Projecting again yields
\begin{equation*}
\Gamma_1=\sum_j(x_j-o_{\widetilde C})\times(\pi_2^X)_*D_j.
\end{equation*}
Its pushforward by $\phi\times\operatorname{id}_X$ is zero, since $\phi_*(x_j-o_{\widetilde C})=0$.
Pushing forward $\Gamma_0+\Gamma_1$ proves the assertion.
\end{proof}

\begin{theorem}\label{thm3.2}
Let $X$ be a smooth projective regular surface with an involution $\sigma$ such that $\mathrm{CH}_0(X/\sigma)=\mathbb Q$.
Assume that there is a degree-one point class $o_X\in\mathrm{CH}^2(X)$ and integral constant-cycle curves $C_\alpha\subseteq X$ whose classes generate $\operatorname{Pic}(X)_{\mathbb Q}$, such that every closed point of every $C_\alpha$ represents this same class $o_X$.
Then $X$ admits motivic $0$-multiplicativity.
\end{theorem}
\begin{proof}
By Lemma \ref{lem3.1} and linearity,
\begin{equation*}
\delta_{X*}A=A\times o_X+o_X\times A
\end{equation*}
for every divisor $A$ on $X$. Intersecting with a divisor on one factor and pushing forward gives
\begin{equation*}
A B=\deg(A B)o_X
\end{equation*}
for all divisors $A,B$.
For an ample divisor $H$, compare $\sigma_*(H^2)$ with $(\sigma_*H)^2$ to obtain $\sigma_*o_X=o_X$.
Since $\mathrm{CH}_0(X/\sigma)=\mathbb Q$, every fixed point of $\sigma$ therefore has class $o_X$.

Let $g:\widetilde X\longrightarrow X$ be the blow-up of the isolated fixed points, and let $\sigma$ also denote the lifted involution.
The quotient $Y:=\widetilde X/\sigma$ is the smooth resolution of $X/\sigma$, and $f:\widetilde X\longrightarrow Y$ is a finite flat double cover, possibly \'etale when the fixed locus is empty.
We use this blow-up construction directly, rather than the scheme-theoretic fiber product with the resolution.
The isomorphism $g_*:\mathrm{CH}_0(\widetilde X)\longrightarrow\mathrm{CH}_0(X)$ shows that every exceptional curve and every strict transform of $C_\alpha$ carries the common class $o_{\widetilde X}:=g^*o_X$.
These curves generate $\operatorname{Pic}(\widetilde X)_{\mathbb Q}$.

The surface $Y$ is regular and has $\mathrm{CH}_0(Y)=\mathbb Q$: resolution of the isolated quotient singularities adds only rational exceptional curves.
Thus $Y$ has trivial Chow groups.
Put $o_Y:=f_*o_{\widetilde X}$. Invariance gives $f^*o_Y=2o_{\widetilde X}$.
For $A\in\mathrm{CH}^1(Y)$, Lemma \ref{lem3.1} gives
\begin{equation*}
\delta_{\widetilde X*}f^*A
=f^*A\times o_{\widetilde X}+o_{\widetilde X}\times f^*A
=\frac12 f_2^*\delta_{Y*}A\in R^*(\widetilde X^2).
\end{equation*}
The last equality follows from the Tate decomposition of the diagonal of the regular surface $Y$.
For codimension-2 cycle classes, which are multiples of $o_Y$, the same ambient membership follows from $f^*o_Y=2o_{\widetilde X}$ and $(\delta_2^{\widetilde X})_*o_{\widetilde X}=o_{\widetilde X}\times o_{\widetilde X}$.
The generalized Lefschetz property holds here: both regular surfaces have vanishing odd cohomology, and $f^*$ is injective in degree two.
Theorem \ref{thm1.1}(iii), in the form of Remark \ref{rem2.2}(ii), therefore makes the natural CK decomposition of $\widetilde X$ multiplicative.
By the criterion of \cite[Proposition~4.1]{FLV21}, $\Gamma^3(\widetilde X,o_{\widetilde X})=0$.
Pushing forward by $g^{\times3}$ gives $\Gamma^3(X,o_X)=0$, so the same criterion supplies an MCK decomposition for $X$.
\end{proof}

\begin{remark}
It is interesting to determine whether the Picard group of every regular surface of geometric genus one is generated by constant-cycle curves. Theorem \ref{thm3.2} requires the additional condition that the generating curves carry one common degree-one class.
\end{remark}

\begin{corollary}
Let $X$ be a smooth projective surface with a relatively minimal elliptic fibration $p:X\longrightarrow\mathbb P^1$ admitting a section.
Then $X$ admits motivic $0$-multiplicativity.
\end{corollary}
\begin{proof}
For a trivial fibration $X\cong E\times\mathbb P^1$, take the product of the standard MCK decompositions.
Otherwise, $X$ is regular and the fibration has a singular fiber.
Fiberwise inversion extends to an involution $\sigma$ of the relatively minimal surface. Its quotient is rational: over $\mathbb C(t)$, a Weierstrass equation $y^2=x^3+a(t)x+b(t)$ has invariant function field $\mathbb C(t,x)$.
Consequently, $\mathrm{CH}_0(X/\sigma)=\mathbb Q$.

Let $O$ be the zero section and let $o_X$ be the class of any point of $O\cong\mathbb P^1$.
Choose a singular fiber $F_0$.
Its integral components are rational curves and form a connected union meeting $O$; their intersection points identify their point classes with $o_X$.
Every section is also a rational curve and meets $F_0$, so all its points have this same class.
The Picard group over $\mathbb Q$ is generated by the zero section, the fiber class, the components of reducible fibers, and sections (the Shioda--Tate description; compare \cite[Remark~2.2]{BV04}).
Replace the fiber class by the weighted sum of components of $F_0$.
Components of every other singular fiber are connected to $O$ in the same way.
Thus the Picard group is generated by integral constant-cycle curves all representing the distinguished class $o_X$.
Therefore the result follows from Theorem \ref{thm3.2}.
\end{proof}

For the next calculations, we use intersection theory for global quotients; refer to \cite{V89}.
If $V=U/G$, with $U$ smooth connected projective and $G$ finite acting faithfully, define its quotient motive by
\begin{equation*}
\mathfrak h(V):=(U,e_G),\qquad e_G:=\frac1{|G|}\sum_{g\in G}\Gamma_g.
\end{equation*}
Rational operational Chow classes on $V$ identify, by cap product with $[V]$, with the corresponding cycle groups; this is the Alexander duality property for quotient singularities.
We use operational pullbacks and these identifications on all quotient powers.
Equivalently, quotient correspondences are interpreted on the invariant summands of $\mathfrak h(U^a)$.
For a quotient map $f:U\longrightarrow U/G$, this formalism gives $f_*f^*=|G|\operatorname{id}$ and $f_2^*\Delta_V=\sum_{g\in G}\Gamma_g$, without requiring $f$ to be flat.

As a consequence of \cite[Remark~4.12(i)]{X25}, we prove the sharp modified-diagonal property for any smooth double EPW sextic, which was predicted by O'Grady for any hyper-Kähler variety in \cite{O'G14}. 
Notably, this property is not a direct application of Theorem \ref{thm1.1}.
We will develop an argument evolving from the approach of \cite[Theorem~4.1]{Voi15}.

\begin{proposition}\label{prop3.5}
Let $Y$ be a general EPW sextic, and let $f: X\rightarrow Y$ be the associated double EPW sextic.
Let $o_X$ be the 0-cycle class represented by any ramification point of $f$.
Then one has the vanishing:
\begin{equation*}
\Gamma^k(X, o_X)=0
\end{equation*}
in $\mathrm{CH}_4(X^k)$ for each integer $k\geq 5$.
\end{proposition}

\begin{proof}
Let $S\subseteq X$ be the ramification locus of $f$, and let $\sigma\colon X\longrightarrow X$ be the covering involution. Thus $S$ is precisely the fixed locus of $\sigma$. By Theorem~1.1 of \cite{Z24}, all closed points of $S$ have the same cycle class in $\mathrm{CH}_0(X)$. We denote this degree-one 0-cycle class by $o_X$.
By \cite[Remark~4.12(i)]{X25}, the Chow motive $\mathfrak h(Y)$ has a natural multiplicative self-dual CK decomposition: 
\begin{equation*}
\pi_{2i}^Y=\frac16 H^{4-i}\times H^i\quad(i=0,1,3,4),\qquad
\pi_{2i+1}^Y=0,\qquad \pi_4^Y=\Delta_Y-\sum_{i\in\{0,1,3,4\}}\pi_{2i}^Y,
\end{equation*}
where $H=c_1(\mathcal O_Y(1))$ and $o_Y=\frac{1}{6}H^4$.
Let $h:=f^*H$. We first prove that
\begin{equation}\label{18}
f^*o_Y=2o_X
\end{equation}

By \cite{F12}, one has
\begin{equation}\label{eq:Ferretti1}
3[S]=15h^2-c_2(X)
\end{equation}
in $\mathrm{CH}^2(X)$,
and
\begin{equation}\label{eq:Ferretti2}
c_2(X)\cdot h=5h^3
\end{equation}
in $\mathrm{CH}^3(X)$.
Multiplying \eqref{eq:Ferretti1} by $h$ and using \eqref{eq:Ferretti2}, we obtain
\begin{equation*}
3[S]\cdot h=10h^3.
\end{equation*}
Hence
\begin{equation}\label{21}
h^3=\frac{3}{10}[S]\cdot h.
\end{equation}
In particular, $h^4$ is represented by a zero-cycle supported on $S$. Since all closed points of $S$ have cycle class $o_X$ and $\deg(h^4)=2\deg(H^4)=12$, we obtain $h^4=12o_X$.
Consequently,
\begin{equation*}
f^*o_Y=\frac{1}{6}h^4=2o_X.
\end{equation*}

Set
\begin{equation}
\Delta_X^+:=(f\times f)^*\Delta_Y=\Delta_X+\Gamma_\sigma
\end{equation}
and
\begin{equation}
\Delta_X^-:=\Delta_X-\Gamma_\sigma.
\end{equation}
Thus
\begin{equation}
2\Delta_X=\Delta_X^++\Delta_X^-.
\end{equation}
By \cite[Lemma~4.3]{Voi15}, one has
\begin{equation}\label{25}
p_{12}^*\Delta_X^-\cdot p_{13}^*\Delta_X^-=p_{12}^*\Delta_X^+\cdot p_{23}^*\Delta_X^-.
\end{equation}
Moreover, since $o_X$ is represented by a fixed point of $\sigma$, one has
\begin{equation}\label{26}
p_2^*o_X\cdot p_{23}^*\Delta_X^-=0.
\end{equation}

We will also use the following consequence of \eqref{21}:
\begin{equation}\label{27}
p_1^*h^3\cdot\Delta_X^-=0
\end{equation}
in $\mathrm{CH}^7(X^2)$,
and
\begin{equation}\label{28}
p_2^*h^3\cdot\Delta_X^-=0
\end{equation}
in $\mathrm{CH}^7(X^2)$.
Indeed, let $j\colon S\hookrightarrow X$ denote the inclusion. By \eqref{21}, there exists a cycle class $\alpha\in\mathrm{CH}^1(S)$ such that $h^3=j_*\alpha$.
Since $\sigma$ restricts to the identity on $S$, one has $\sigma\circ j=j$.
Consequently, the restrictions to $S$ of the diagonal embedding and of the graph embedding of $\sigma$ coincide. It follows that
\begin{equation}
p_1^*h^3\cdot\Delta_X
=p_1^*h^3\cdot\Gamma_\sigma,
\end{equation}
which proves \eqref{27}. The proof of \eqref{28} is identical.

Put
\begin{equation}
\widetilde{\Delta}_Y:=\Delta_Y-p_2^*o_Y.
\end{equation}
By \eqref{18}, one has
\begin{equation}
(f\times f)^*\widetilde{\Delta}_Y
=\Delta_X^+-2p_2^*o_X.
\end{equation}
Combining this equality with \eqref{25} and \eqref{26}, we obtain
\begin{equation}\label{eq:pairing}
p_{12}^*\Delta_X^-\cdot p_{13}^*\Delta_X^-
=p_{12}^*(f\times f)^*\widetilde{\Delta}_Y
\cdot p_{23}^*\Delta_X^-.
\end{equation}

By \cite{X25}, the fourfold $Y$ satisfies
\begin{equation}\label{33}
\Gamma^3(Y,o_Y)
\in\operatorname{Span}_{\mathbb{Q}}
\langle H_1^2H_2^3H_3^3, H_1^3H_2^2H_3^3, H_1^3H_2^3H_3^2\rangle,
\end{equation}
where $H_i$ denotes the pullback of $H$ from the $i$-th factor. Likewise,
\begin{equation}\label{34}
\Gamma^4(Y,o_Y)
\in
\mathbb{Q}\cdot H_1^3H_2^3H_3^3H_4^3.
\end{equation}
Moreover,
\begin{equation}\label{35}
\Gamma^k(Y,o_Y)=0
\end{equation}
for every integer $k\geq5$.

For every integer $k\geq3$, define
\begin{equation}\label{eq:def-Gamma1r}
\Gamma^{1,k}(Y,o_Y)
:=
\prod_{i=1}^{k}q_{0i}^*\widetilde{\Delta}_Y
\in
\mathrm{CH}_4(Y^{k+1}),
\end{equation}
where $q_{0i}\colon Y^{k+1}\longrightarrow Y^2$ denotes the projection onto the $0$-th and $i$-th factors. The standard relation between $\Gamma^{1,k}(Y,o_Y)$ and the usual modified diagonal class is
\begin{equation}\label{eq:Gamma1r-relation}
\Gamma^{1,k}(Y,o_Y)
=\Gamma^{k+1}(Y,o_Y)
+q_0^*o_Y\cdot q_{1\cdots k}^*\Gamma^k(Y,o_Y).
\end{equation}
In particular,
\begin{equation}\label{38}
\Gamma^{1,3}(Y,o_Y)
=\Gamma^4(Y,o_Y)
+q_0^*o_Y\cdot q_{123}^*\Gamma^3(Y,o_Y).
\end{equation}
Using \eqref{35}, we obtain
\begin{equation}\label{39}
\Gamma^{1,4}(Y,o_Y)
=q_0^*o_Y\cdot q_{1234}^*\Gamma^4(Y,o_Y)
\end{equation}
and
\begin{equation}\label{40}
\Gamma^{1,5}(Y,o_Y)=0.
\end{equation}

We now turn to $X$. Define
\begin{equation}\label{41}
\Gamma^{1,5}(X,o_X):=\prod_{i=1}^{5}\left(p_{0i}^*\Delta_X-p_i^*o_X\right).
\end{equation}
in $\mathrm{CH}_4(X^6)$.
Let $p'_{1\cdots5}\colon X^6\longrightarrow X^5$ be the projection onto the last five factors. 
Then
\begin{equation}\label{42}
\Gamma^5(X,o_X)
=(p'_{1\cdots5})_*\Gamma^{1,5}(X,o_X).
\end{equation}

By \eqref{18}, one has
\begin{equation}\label{43}
2\left(\Delta_X-p_2^*o_X\right)
=(f\times f)^*\widetilde{\Delta}_Y+\Delta_X^-.
\end{equation}
Consequently,
\begin{equation}\label{44}
2^5\Gamma^{1,5}(X,o_X)
=\prod_{i=1}^{5}
\left(
p_{0i}^*(f\times f)^*\widetilde{\Delta}_Y
+
p_{0i}^*\Delta_X^-
\right).
\end{equation}

For convenience, set
$F_{ij}:=p_{ij}^*(f\times f)^*\widetilde{\Delta}_Y$ and $E_{ij}:=p_{ij}^*\Delta_X^-$.
We further write $F_i:=F_{0i}$ and $E_i:=E_{0i}$.
Equation \eqref{eq:pairing} gives
\begin{equation}\label{45}
E_iE_j=F_iE_{ij}
\end{equation}
for $i\neq j$.
Moreover, equations \eqref{26}, \eqref{27}, and \eqref{28} imply
\begin{equation}\label{46}
p_i^*o_X\cdot E_{ij}=0
\end{equation}
and
\begin{equation}\label{47}
p_i^*h^3\cdot E_{ij}=0.
\end{equation}

Upon expanding the right-hand side of \eqref{44}, every summand containing at least one $F_i$ is, up to a permutation of the last five factors, of the form
\begin{equation}\label{eq:def-Tk}
T_k:=
F_1\cdots F_kE_{k+1}\cdots E_5
\end{equation}
for some $1\leq k\leq5$. We claim that
\begin{equation}\label{eq:Tk-vanishing}
T_k=0
\end{equation}
for every integer $1\leq k\leq 5$.

For $k=1$, equation \eqref{45} gives
$E_2E_3=F_2E_{23}$ and $E_4E_5=F_4E_{45}$.
Therefore
\begin{equation*}
T_1=F_1E_2E_3E_4E_5=F_1F_2F_4E_{23}E_{45}.
\end{equation*}
The class $F_1F_2F_4$ is the pullback of $\Gamma^{1,3}(Y,o_Y)$ from the factors $0,1,2,4$. 
By \eqref{38}, it is the sum of the pullback of $\Gamma^4(Y,o_Y)$ and a cycle class which, 
up to a non-zero rational scalar, is of the form
\begin{equation}
p_0^*o_X\cdot
(f_3\circ p_{124})^*\Gamma^3(Y,o_Y).
\end{equation}
By \eqref{34}, the first part is a rational multiple of a cycle class containing $p_2^*h^3\cdot p_4^*h^3$.
Hence its product with $E_{23}E_{45}$ vanishes by \eqref{47}.
On the other hand, by \eqref{33}, every monomial occurring in the second part has, on the factors $1,2,4$, exponents given by a permutation of $(2,3,3)$.
Thus at least one of the factors $2$ and $4$ carries a factor $h^3$. Therefore the corresponding term is annihilated by $E_{23}$ or $E_{45}$, again by \eqref{47}. Consequently, $T_1=0$.

For $k=2$, equation \eqref{45} gives $E_3E_4=F_3E_{34}$.
Hence
\begin{equation*}
T_2=F_1F_2E_3E_4E_5=F_1F_2F_3E_{34}E_5.
\end{equation*}
The product $F_1F_2F_3$ is the pullback of $\Gamma^{1,3}(Y,o_Y)$. By \eqref{34}, its $\Gamma^4(Y,o_Y)$-part contains a factor $p_3^*h^3$ and is therefore annihilated by $E_{34}$, by \eqref{47}. Its $q_0^*o_Y\cdot\Gamma^3(Y,o_Y)$-part contains a factor $p_0^*o_X$ and is therefore annihilated by
$E_5=E_{05}$, in view of \eqref{46}. Thus $T_2=0$.

For $k=3$, equation \eqref{45} gives $E_4E_5=F_4E_{45}$.
Hence
\begin{equation*}
T_3=F_1F_2F_3E_4E_5=F_1F_2F_3F_4E_{45}.
\end{equation*}
By \eqref{39}, the product $F_1F_2F_3F_4$ is, up to a non-zero rational scalar, 
the pullback of $q_0^*o_Y\cdot q_{1234}^*\Gamma^4(Y,o_Y)$.
By \eqref{34}, this cycle class contains a factor $p_4^*h^3$. Hence it is annihilated by $E_{45}$, by \eqref{47}. 
Therefore $T_3=0$.

For $k=4$, one has
\begin{equation*}
T_4=F_1F_2F_3F_4E_5.
\end{equation*}
By \eqref{39}, the product $F_1F_2F_3F_4$ contains, up to a non-zero rational scalar, a factor $p_0^*o_X\cdot(f_4\circ p_{1234})^*\Gamma^4(Y,o_Y)$. 
Since $E_5=E_{05}$, equation \eqref{46} gives $T_4=0$.

Finally,
\begin{equation}\label{eq:T5}
T_5=F_1F_2F_3F_4F_5.
\end{equation}
This is the pullback of $\Gamma^{1,5}(Y,o_Y)$. Therefore \eqref{40} gives
$T_5=0$. This proves \eqref{eq:Tk-vanishing}.

It follows from \eqref{44} that
\begin{equation}\label{eq:pure-minus}
2^5\Gamma^{1,5}(X,o_X)=\prod_{i=1}^{5}p_{0i}^*\Delta_X^-.
\end{equation}

Let $\sigma':=(\sigma,\operatorname{id}_{X^5}): X^6\longrightarrow X^6$.
For every $i$, one has $(\sigma')^*p_{0i}^*\Delta_X^-=-p_{0i}^*\Delta_X^-$.
Since there are five factors, it follows that
\begin{equation}\label{eq:anti-invariance}
(\sigma')^*
\left(
\prod_{i=1}^{5}p_{0i}^*\Delta_X^-
\right)
=-\prod_{i=1}^{5}p_{0i}^*\Delta_X^-.
\end{equation}
Since $\sigma'$ is an involution, one also has
\begin{equation}\label{54}
(\sigma')_*
\left(
\prod_{i=1}^{5}p_{0i}^*\Delta_X^-
\right)
=-\prod_{i=1}^{5}p_{0i}^*\Delta_X^-.
\end{equation}
Moreover, $p'_{1\cdots5}\circ\sigma'=p'_{1\cdots5}$.
Using the anti-invariance \eqref{54}, we then obtain
\begin{equation*}
(p'_{1\cdots5})_*
\left(
\prod_{i=1}^{5}p_{0i}^*\Delta_X^-
\right)
=(p'_{1\cdots5})_*(\sigma')_*
\left(
\prod_{i=1}^{5}p_{0i}^*\Delta_X^-
\right)=-(p'_{1\cdots5})_*\left(\prod_{i=1}^{5}p_{0i}^*\Delta_X^-\right).
\end{equation*}
It follows that
\begin{equation}\label{eq:push-pure-minus-zero}
(p'_{1\cdots5})_*\left(\prod_{i=1}^{5}p_{0i}^*\Delta_X^-\right)=0.
\end{equation}
Combining \eqref{42}, \eqref{eq:pure-minus}, and \eqref{eq:push-pure-minus-zero}, we obtain
\begin{equation*}
2^5\Gamma^5(X,o_X)=0.
\end{equation*}
Since the Chow groups are with rational coefficients, then
\begin{equation*}
\Gamma^5(X,o_X)=0
\end{equation*}
in $\mathrm{CH}_4(X^5)$. 
By \cite[Proposition~2.4]{O'G14}, this vanishing propagates to higher modified diagonals. Therefore
\begin{equation*}
\Gamma^k(X,o_X)=0
\end{equation*}
in $\mathrm{CH}_4(X^k)$ for every integer $k\geq 5$.
This completes the proof.
\end{proof}

Finally, we promote a result in \cite{BL24} to Chow, which was only proved modulo algebraic equivalence.
\begin{theorem}\label{thm:GM-sixfold-MCK}
Let $X$ be a smooth Gushel--Mukai sixfold. Then the natural self-dual CK decomposition of Theorem \ref{thm1.1} is multiplicative.
\end{theorem}
\begin{proof} 
For a smooth projective variety $V$ and $\alpha\in\mathrm{CH}^r(V)$, denote
\begin{equation*}
L_{\alpha}:=\Delta_V\cdot p_2^*\alpha\in\mathrm{CH}^{\dim V+r}(V^2)
\end{equation*}
for multiplication by $\alpha$, regarded as a correspondence. If $\tau_V$ is an involution of $V$, put $\pi_V^-:=\frac12(\Delta_V-\Gamma_{\tau_V})$ and denote its anti-invariant cohomology by $H^*(V)^-$.

\emph{Step 1: reduce to one weighted diagonal.}
Let $G=\operatorname{Gr}(2,V_5)$, let $f:X\longrightarrow G$ be the Gushel double cover, and let $\sigma$ be its involution. Put
\begin{equation*}
H:=\sigma_1,\qquad C:=\sigma_{1,1},\qquad h:=f^*H,\qquad c:=f^*C.
\end{equation*}
The reduced ramification divisor has class $h$ and is fixed pointwise by $\sigma$. The cycle equality in proof of Theorem \ref{thm1.1} therefore gives
\begin{equation}\label{eq:GM-h-supported}
L_{h f^*\beta}\circ\pi_X^-=0
\qquad\bigl(\beta\in\mathrm{CH}^*(G)\bigr).
\end{equation}
The Chow ring of $G$ has the presentation:
\begin{equation*}
\mathrm{CH}^*(G)=\frac{\mathbb Q[H,C]}{
(H^4-3H^2C+C^2,\ H^5-4H^3C+3HC^2)}.
\end{equation*}
Consequently, $\mathrm{CH}^*(G)/(H)=\mathbb Q[\overline C]/(\overline C^{\,2})$.
For $\alpha\in\mathrm{CH}^*(G)$, the double-cover identity gives
\begin{equation}\label{eq:GM-weighted-diagonal}
\begin{split}
T_\alpha&:=\delta_{X*}f^*\alpha-\frac12(f\times f)^*\delta_{G*}\alpha\\
&=L_{f^*\alpha}\circ\pi_X^-.
\end{split}
\end{equation}
By \eqref{eq:GM-h-supported}, every $T_\alpha$ with $\alpha$ of positive codimension is a multiple of the cycle class
\begin{equation}\label{eq:GM-Theta}
\Theta_X:=T_C=L_c\circ\pi_X^-\in\mathrm{CH}^8(X^2).
\end{equation}
Thus Theorem \ref{thm1.1}(iii) reduces the assertion to $\Theta_X=0$.

\emph{Step 2: use the Franchetta property in lower codimensions.}
Let $\mathcal X\longrightarrow B$ be the family of smooth quadric sections of the cone over $G$, as in \cite[Definition~4.3]{BL24}. For a fiber, set $R^*(X):=f^*\mathrm{CH}^*(G)$ and let $R^{(2),*}$ be the image of $R^*(X)\otimes R^*(X)$ in $\mathrm{CH}^*(X^2)$. Write $\operatorname{GDCH}_B^*$ for the image of the Chow groups of the total relative family on a fiber.
The stratified projective-bundle calculation of \cite[Equation~(4.7)]{BL24} gives
\begin{equation}\label{eq:GM-GD-square}
\operatorname{GDCH}_B^*(X^2)=\langle R^{(2),*},\Delta_X\rangle.
\end{equation}
The cycle class map is injective on $R^{(2),*}$, by the Tate motive of $G$ and $f_*f^*=2\operatorname{id}$.
In codimension $5$, the right-hand side of \eqref{eq:GM-GD-square} is ambient. In codimension $6$, its additional diagonal class is independent of the ambient classes in cohomology, since it acts as the identity on $H^6(X)^-$; this is also \cite[Proposition~4.6]{BL24}. In codimension $7$, the only additional generator is $\Delta_X\cdot p_2^*h$, and
\begin{equation*}
\Delta_X\cdot p_2^*h
=\Gamma_\sigma\cdot p_2^*h
=\frac12(f\times f)^*(\Delta_G\cdot p_2^*H).
\end{equation*}
Hence,
\begin{equation}\label{eq:GM-low-Franchetta}
\operatorname{GDCH}_B^r(X^2)\hookrightarrow H^{2r}(X^2,\mathbb Q)
\qquad(r=5,6,7).
\end{equation}

Let $\mathcal U_X:=f^*\mathcal U_G$, where $\mathcal U_G$ is the tautological rank-two subbundle and let
\begin{equation*}
q:E:=\mathbb P_X(\mathcal U_X)\longrightarrow X
\end{equation*}
parametrizes lines in its fibers. The first quadric fibration is $\rho:E\longrightarrow\mathbb P(V_5)$; put $u:=\rho^*c_1\mathcal O_{\mathbb P(V_5)}(1)=c_1\mathcal O_E(1)$.
The projective-bundle formula for the relative square has basis $1,u_1,u_2,u_1u_2$. Its coefficients in codimension $7$ have codimensions $7,6,6,5$ on $X^2$. Therefore \eqref{eq:GM-low-Franchetta} implies
\begin{equation}\label{eq:GM-E-Franchetta}
\operatorname{GDCH}_B^7(E^2)\hookrightarrow H^{14}(E^2,\mathbb Q).
\end{equation}
These assertions remain valid after replacing $B$ by a Zariski open subset. 

\emph{Step 3: establish an incidence isomorphism in Chow.}
We now work on a dense open subset of $B$, where the associated double EPW sextic $Y=\widetilde Y_A$ is smooth and the hyperplane section
\begin{equation*}
j:D:=\widetilde Y_A\times_{\mathbb P(V_6)}\mathbb P(V_5)\hookrightarrow Y
\end{equation*}
is smooth. We impose the generality condition (10) of \cite{DK19}. Put $\ell=c_1\mathcal O_Y(1)$ and $\ell_D=j^*\ell$.
The constructions below are made over the original quadric parameter space. In particular, $D$ is intrinsically the Stein factor of the scheme of $\sigma$-planes over the discriminant of $\rho$. Over the algebraic closure of the function field of $B$, identify its hyperplane section with this intrinsic $D$, as in \cite[Theorem~4.3(a)]{DK19}. 
Clearly, the cover $Y$ together with $j$ descends to the original function field and spreads over a Zariski open subset of $B$.

By \cite[Theorem~4.3(a) and Corollary~5.13]{DK19}, the scheme $F=F_2^\sigma(X)$ is smooth and admits an \'etale locally trivial $\mathbb P^1$-fibration $r:F\longrightarrow D$. Let $p:\mathcal P\longrightarrow F$ be the universal plane and let $v:\mathcal P\longrightarrow E$ be its natural lift. Put
\begin{equation*}
\xi:=-\frac12c_1(\omega_{F/D})\in\mathrm{CH}^1(F),\qquad r_*\xi=[D],
\end{equation*}
so no choice of a section of $r$ is needed, and define
\begin{equation}\label{eq:GM-incidence}
\Gamma:=\pi_E^-\circ(rp,v)_*(p^*\xi)\circ\pi_D^-
\in\mathrm{CH}^5(D\times E).
\end{equation}
Here $\sigma$ lifts to $E$ and to $F$; on the two rulings of a general singular quadric fiber it induces the nontrivial deck involution of $D$. The polarization classes and $\xi$ are invariant under these involutions. Thus, all the projectors in \eqref{eq:GM-incidence} preserve the support over $\mathbb P(V_5)$.

The support of $\Gamma$ lies in $D\times_{\mathbb P(V_5)}E$. Both projections to $\mathbb P(V_5)$ pull back the same hyperplane class on this support, so
\begin{equation}\label{eq:GM-incidence-linear}
L_u\circ\Gamma=\Gamma\circ L_{\ell_D},\qquad
{}^t\Gamma\circ L_u=L_{\ell_D}\circ{}^t\Gamma.
\end{equation}
Moreover, ${}^t\Gamma\circ\Gamma\in\mathrm{CH}^3(D^2)$ is supported on $D\times_{\mathbb P(V_5)}D$. This finite fiber product has only the diagonal and the involution graph as its three-dimensional components. 
So, a 3-dimensional cycle class on it is a $\mathbb{Q}$-linear combination of these components. Since the composition is projected to the anti-invariant summand on both sides, localization gives the Chow identity
\begin{equation}\label{eq:GM-incidence-left}
{}^t\Gamma\circ\Gamma=\lambda\pi_D^-
\qquad\text{for some }\lambda\in\mathbb Q.
\end{equation}

We verify that $\lambda\ne0$. The projective-bundle formula and the generalized Lefschetz theorem give
\begin{equation*}
H^6(E)^-=q^*H^6(X)^-,\qquad
H^8(E)^-=u\,q^*H^6(X)^-,
\end{equation*}
with no other anti-invariant cohomology on $E$.
The plane Abel--Jacobi isomorphism of \cite[Theorem~5.19]{DK19}, together with $r_*\xi=1$, identifies
\begin{equation}\label{eq:GM-incidence-cohomology}
({}^t\Gamma)_*:H^6(E)^-\xrightarrow{\ \sim\ }H^2(D)^-.
\end{equation}
Indeed, for $x\in H^6(X)^-$ the plane Abel--Jacobi image $p_*v^*q^*x$ is $r^*\beta_D(x)$, where $\beta_D$ is the restriction of that isomorphism to $D$; multiplying by $\xi$ and pushing along $r$ gives $\beta_D(x)$.
Using \eqref{eq:GM-incidence-linear}, the projective-bundle formula, and hard Lefschetz on the smooth threefold $D$, the same correspondence is an isomorphism $H^8(E)^-\longrightarrow H^4(D)^-$.
Poincar\'e duality then makes $\Gamma_*:H^2(D)^-\longrightarrow H^6(E)^-$ and $\Gamma_*:H^4(D)^-\longrightarrow H^8(E)^-$ isomorphisms. Thus the composition in \eqref{eq:GM-incidence-left} is nonzero on $H^2(D)^-$ and $\lambda\ne0$.
Equation \eqref{eq:GM-incidence-left} also forces $H^k(D)^-=0$ for $k\ne2,4$, since $\Gamma_*$ then maps to a zero anti-invariant group of $E$.
It follows that $\Gamma\circ{}^t\Gamma=\lambda\pi_E^-$ in cohomology. This difference is a generically defined codimension-7 cycle class on $E^2$. By \eqref{eq:GM-E-Franchetta},
\begin{equation}\label{eq:GM-incidence-right}
\Gamma\circ{}^t\Gamma=\lambda\pi_E^-
\quad\text{in }\mathrm{CH}^7(E^2).
\end{equation}

\emph{Step 4: split the hyperplane inclusion on the anti-invariant motive.}
For a smooth double EPW sextic, $H^2(Y)^+=\mathbb Q\ell$ and $H^4(Y)=\operatorname{Sym}^2H^2(Y)$. Consequently,
\begin{equation*}
H^4(Y)^-=\ell H^2(Y)^-,\qquad H^6(Y)^-=\ell^2H^2(Y)^-.
\end{equation*}
The Lefschetz standard conjecture for varieties of $K3^{[2]}$-type holds by \cite[Theorem~1.1]{CM13}. Choose $\Psi\in\mathrm{CH}^2(Y^2)$, projected to the anti-invariant summands, whose cohomological action inverts $\ell^2:H^2(Y)^-\longrightarrow H^6(Y)^-$, and put
\begin{equation}\label{eq:GM-Lefschetz-inverse}
\Lambda:=\Psi\circ L_\ell+L_\ell\circ\Psi\in\mathrm{CH}^3(Y^2).
\end{equation}
On $H^4(Y)^-$ and $H^6(Y)^-$, respectively, $\Lambda_*$ is the inverse of multiplication by $\ell$. The other term in each degree is zero, since $H^0(Y)^-=H^8(Y)^-=0$.

We make the spread of this auxiliary correspondence explicit. On a very general fiber, the algebraic inverse supplied by \cite{CM13} belongs to a component of the relative parameter space of cycles dominating $B$; its prescribed cohomological action is the horizontal inverse to $\ell^2$.
Taking a generically finite multisection and shrinking $B$ gives a finite \'etale cover $B'\longrightarrow B$ over which the cycle spreads. If its degree is $N$, replace the spread by
\begin{equation*}
\frac1N(\mathcal Y_{B'}\times_{B'}\mathcal Y_{B'}
\longrightarrow\mathcal Y\times_B\mathcal Y)_*\Psi_{B'}.
\end{equation*}
Its cohomological action is still the prescribed inverse. Project it again by $\pi_Y^-$, if necessary, and use it in \eqref{eq:GM-Lefschetz-inverse}. Thus $\Psi$, $\Lambda$, and the compositions below are defined over an open subset of the original $B$.

Define
\begin{equation}\label{eq:GM-hyperplane-J}
J:=\pi_D^-\circ{}^t\Gamma_j\circ\Lambda\circ\Gamma_j\circ\pi_D^-
\in\mathrm{CH}^3(D^2).
\end{equation}
Weak Lefschetz identifies $H^2(D)^-$ with $j^*H^2(Y)^-$. For $a\in H^2(Y)^-$, the projection formula gives
\begin{equation*}
j_*j^*a=\ell a,\qquad
j_*(\ell_D j^*a)=\ell^2a.
\end{equation*}
Together with \eqref{eq:GM-Lefschetz-inverse}, these identities give $J=\pi_D^-$ in cohomology in degrees $2$ and $4$, hence in every degree by the vanishing established above.
The generically defined cycle
\begin{equation*}
\Gamma\circ(J-\pi_D^-)\circ{}^t\Gamma
\in\mathrm{CH}^7(E^2)
\end{equation*}
is therefore homologically trivial, and \eqref{eq:GM-E-Franchetta} makes it zero in Chow. Composing on the left with ${}^t\Gamma$ and on the right with $\Gamma$, and using \eqref{eq:GM-incidence-left}, gives $\lambda^2(J-\pi_D^-)=0$. Hence $J=\pi_D^-$ in Chow.
In particular, the correspondence
\begin{equation*}
\kappa:=\pi_D^-\circ{}^t\Gamma_j\circ\Lambda\circ\pi_Y^-
\in\mathrm{CH}^3(Y\times D)
\end{equation*}
satisfies
\begin{equation}\label{eq:GM-hyperplane-split}
\kappa\circ\Gamma_j\circ\pi_D^-=\pi_D^-.
\end{equation}
This is a left inverse to the Gysin map on the anti-invariant Chow motive.

\emph{Step 5: transfer to the quadric fibration.}
Let $S\subseteq Y$ be the fixed surface of its covering involution. 
The cycle equalities of \cite[Lemma~4.1 and Proposition~4.5]{F12} give
\begin{equation*}
3[S]=15\ell^2-c_2(Y),\qquad c_2(Y)\ell=5\ell^3,
\end{equation*}
and thus $[S]\ell=\frac{10}{3}\ell^3$ in Chow. 
Since the involution is the identity on $S$, the diagonal and involution graph have identical restrictions to this supporting cycle. Therefore,
\begin{equation}\label{eq:GM-EPW-fixed}
L_{\ell^3}\circ\pi_Y^-=0\quad\text{in }\mathrm{CH}^7(Y^2).
\end{equation}
Put $w_D:=L_{\ell_D^2}\circ\pi_D^-\in\mathrm{CH}^5(D^2)$. It is anti-invariant in both factors, and the projection formula gives
\begin{equation*}
(j\times j)_*w_D=L_{\ell^3}\circ\pi_Y^-=0.
\end{equation*}
Applying $\kappa\otimes\kappa$ and \eqref{eq:GM-hyperplane-split} yields
\begin{equation}\label{eq:GM-D-vanishing}
L_{\ell_D^2}\circ\pi_D^-=0\quad\text{in }\mathrm{CH}^5(D^2).
\end{equation}
Now \eqref{eq:GM-incidence-linear}, \eqref{eq:GM-incidence-right}, and \eqref{eq:GM-D-vanishing} imply
\begin{equation}\label{eq:GM-E-u2}
L_{u^2}\circ\pi_E^-
=\frac1\lambda\Gamma\circ
(L_{\ell_D^2}\circ\pi_D^-)\circ{}^t\Gamma=0.
\end{equation}

The tautological sequence on $E$, with subbundle $\mathcal O_E(-1)\subset q^*\mathcal U_X$, gives
\begin{equation}\label{eq:GM-projective-relation}
q^*c=u(q^*h-u).
\end{equation}
The inverse image under $q$ of the ramification divisor of $f$ is fixed pointwise by the lifted involution and has class $q^*h$. Hence $L_{q^*h}\circ\pi_E^-=0$ in Chow. Combining this with \eqref{eq:GM-E-u2} and \eqref{eq:GM-projective-relation} gives
\begin{equation*}
L_{q^*c}\circ\pi_E^-=0.
\end{equation*}
Finally, the projective-bundle equality $q_*(u\cdot q^*(-))=\operatorname{id}$ holds as a Chow-correspondence identity, and $q^*$ intertwines the involution projectors. Composing the preceding vanishing with $q^*$ and its left inverse $q_*(u\cdot-)$ therefore gives
\begin{equation}\label{eq:GM-Theta-zero}
\Theta_X=L_c\circ\pi_X^-=0\quad\text{in }\mathrm{CH}^8(X^2).
\end{equation}

The argument proves \eqref{eq:GM-Theta-zero} on the geometric generic fiber of the smooth quadric family. 
By specialization, \eqref{eq:GM-Theta-zero} holds for any smooth GM sixfold.
Together with reduction of Step 1 and Theorem \ref{thm1.1}(iii), this proves the multiplicativity of the natural CK decomposition.
\end{proof}

\begin{remark}\label{rem:GM-Franchetta-surfaces}
Our proof uses the Franchetta property for $X^2$ only in codimensions $5,6,7$, but not in codimension $8$.
In fact, \eqref{eq:GM-GD-square} and \eqref{eq:GM-weighted-diagonal} give
\begin{equation*}
\operatorname{GDCH}_B^8(X^2)=R^{(2),8}+\mathbb Q\Theta_X.
\end{equation*}
Thus \eqref{eq:GM-Theta-zero} proves precisely the missing codimension-8 Franchetta property. The new ingredient is the Chow-level splitting \eqref{eq:GM-hyperplane-split}, which transfers the fixed-surface relation \eqref{eq:GM-EPW-fixed}. 

The cancellation of the six surface classes can also be made explicitly. For a general hyperplane $V_4\subset V_5$, put $Q=\operatorname{Gr}(2,V_4)$ and let $i:Z=f^{-1}(Q)\hookrightarrow X$. Since $[Q]=C$, one has $i_*[Z]=c$. The smooth fourfold $Z$ is an intersection of two quadrics in $\mathbb P^6$.
By \cite[Corollary~2.1]{BT16}, we have a decomposition
\begin{equation*}
\mathfrak h(Z)\cong\mathbf1\oplus\mathbb L\oplus
(\mathbb L^2)^{\oplus8}\oplus\mathbb L^3\oplus\mathbb L^4.
\end{equation*}
The invariant motive is $\mathfrak h(Q)$, whose middle cohomology has rank two.
Hence, we have $\mathfrak h(Z)^-\cong(\mathbb L^2)^{\oplus 6}$.
Choose a basis $\alpha_1,\ldots, \alpha_6$ of $\mathrm{CH}^2(Z)^-$, put $c_{ij}=\deg(\alpha_i\alpha_j)$, and let $(c^{ij})$ be the inverse intersection matrix. The previous decomposition gives
\begin{equation*}
\pi_Z^-=\sum_{i,j=1}^6c^{ij}\, \alpha_i\times\alpha_j.
\end{equation*}
Denoting $\beta_i=i_*\alpha_i\in\mathrm{CH}^4(X)_{\mathrm{hom}}^-$, then we deduce that
\begin{equation}\label{eq:GM-six-surfaces-zero}
\sum_{i,j=1}^6c^{ij}\, \beta_i\times\beta_j
=(i\times i)_*\pi_Z^-=\Theta_X=0
\quad\text{in }\mathrm{CH}^8(X^2).
\end{equation}
\end{remark}
\bigskip
\bigskip
\bigskip

{\noindent\textbf{Funding}}
\medskip

This work was supported by the Fundamental Research Funds of Shandong University
under Grant No.~11140075614046, the Future Program for Young Scholars of Shandong University under Grant No.~11140089964236, and the National Natural Science Foundation of China under Grant Nos.~11601276 and 11771426.
\bigskip
\bigskip

{\noindent\textbf{Acknowledgments}}
\medskip

The author would like to express sincere gratitude to Professors Baohua Fu and Kejian Xu for their constant encouragement and support.
\bigskip
\bigskip

{\noindent\textbf{Declaration of AI Use.}} The author used ChatGPT (OpenAI) to assist with language editing and preliminary mathematical discussions. All mathematical statements and proofs were independently verified by the author, who takes full responsibility for the manuscript.

\vspace{1cm}

{\small

\noindent
{\bf Ze Xu}\\ School of Mathematics, Shandong University, 27 South Shanda Road, Jinan, Shandong 250100, P. R. China\\
{\bf Email: xuze@sdu.edu.cn}

}
\end{document}